\documentclass[11pt]{article}
\usepackage{amsmath,amssymb}
\usepackage{epsfig}
\usepackage{amsfonts}
\usepackage{amsmath,amssymb}
\usepackage{epsfig}
\def\dref#1{(\ref{#1})}

\begin{document}

\vspace*{1\baselineskip} \begin{center}\Large\bf Disconnected
synchronized regions of complex dynamical networks \footnote {\small
This work is supported by the National
 Science Foundation of
 China under grants 60674093, 60334030. }
 \end{center}
\vspace*{1\baselineskip}

\centerline {Zhisheng Duan,  ~Guanrong Chen ~and ~Lin Huang}

\vspace*{0.5\baselineskip}
\begin{center}
 State Key Laboratory for Turbulence and Complex Systems,
Department of Mechanics and Aerospace Engineering, College of
Engineering, Peking University, Beijing 100871, P. R. China
 \\  Email:
       duanzs@pku.edu.cn, eegchen@cityu.edu.hk,
       hl35hj75@pku.edu.cn\\
       {\it Tel and Fax}: (8610)62765037
\end{center}
\vskip 0.5cm

{\bf Abstract.} \,\, This paper addresses the synchronized region
problem, which is reduced to a matrix stability problem, for complex
dynamical networks. For any natural number $n$, the existence of a
network which has $n$ disconnected synchronized regions is
theoretically demonstrated. This shows the complexity in network
synchronization. Convexity characteristic of stability for matrix
pencils is further discussed. Smooth and generalized smooth Chua's
circuit networks are finally discussed as examples for illustration.


\section{ Introduction and problem formulation}

 The subject of network synchronization has recently attracted increasing attention
 from various fields
(see \cite{bar02,bel05, boc06, cur97, koc05, luw04, lu04, pec98,
sor07, wang06, wang02, wat98, wu02, zhao06} and references therein).
 Of particular interest is how the synchronization ability depends on
 various structural parameters of the network, such as average
  distance,
 clustering coefficient, coupling strength, degree distribution and weight distribution.
Some important results have been established for such problems by
introducing the notions of master stability function and
synchronized region \cite{bar02,koc05,lu04,mot05,pec98,zhou06}. It
is natural to expect strong synchronization ability at small cost
\cite{nis07}. In fact, a key factor influencing the synchronization
ability is the characterization of the network synchronized region,
as studied in \cite{koc05, liu07, pec98}. Obviously, the larger the
synchronized region, the easier the synchronization. Some examples
for the existence of two and three disconnected synchronized regions
are demonstrated in \cite{liu07}. This paper attempts to explore the
existence of multiple disconnected synchronized regions for various
complex dynamical networks.

Consider a dynamical network consisting of $N$ coupled identical
nodes, with each node being an $n$-dimensional dynamical system,
described by
\begin{equation}\label{n1}
\dot{x}_i=f(x_i)+c\sum_{j=1}^Na_{ij}H(x_j),\;i=1,2,\cdots,N,
\end{equation}
where $x_i=(x_{i1},x_{i2},\cdots,x_{in})\in \mathbb{R}^n$ is the
state vector of node $i$, $f(\cdot):\mathbb{R}^n\rightarrow
\mathbb{R}^n$
is a smooth 
vector-valued function,  constant $c>0$ represents the
\textit{coupling strength},  $H(\cdot):\mathbb{R}^n\rightarrow
\mathbb{R}^n$ is called the
\textit{inner linking function}, 
and $A=(a_{ij})_{N\times N}$ is called the \textit{outer coupling
matrix}, which represents the coupling configuration of the entire
network. Generally,  $A$ is an irreducible matrix, and if the
entries of $A$ satisfy
$$
a_{ii}=-\sum_{j=1,j\neq i}^Na_{ij},\;i=1,2,\cdots,N,
$$
then network (\ref{n1}) is called a diffusively coupled network. In
this case,  zero is an eigenvalue of $A$ with  multiplicity $1$ and
all the other eigenvalues of $A$ are strictly  negative, which are
denoted by
 \begin{equation}\label{f1}
0=\lambda_1>\lambda_2\geq\lambda_3\geq\cdots\geq\lambda_N.
\end{equation}

The dynamical network \dref{n1} is said to achieve (asymptotical)
synchronization if
\begin{equation}\label{f2}
x_1(t)\rightarrow x_2(t)\rightarrow\cdots\rightarrow
x_N(t)\rightarrow s(t),\;\textrm{as}\; t\rightarrow \infty,
\end{equation}
where, because of the diffusive coupling configuration,  the
\textit{synchronous state} $s(t)\in \mathbb{R}^n$ is a solution of
an individual node, i.e., $\dot{s}(t)=f(s(t))$. Here, $s(t)$ can be
an equilibrium point, a periodic orbit, or even a chaotic orbit.

As shown in \cite{lu04,pec98}, the stability of the synchronized
solution $x_1(t)=x_2(t)=\cdots=x_N(t)=s(t)$
can be determined by analyzing the following equation, known as the
\textit{master stability equation}:
\begin{equation}\label{f3}
\dot{\omega}=[Df(s(t))+\alpha DH(s(t))]\omega,
\end{equation}
where $\alpha\in \mathbb{R}$, and $Df(s(t))$ and $DH(s(t))$ are the
Jacobian matrices of functions $f$ and $H$ at $s(t)$, respectively.

The largest Lyapunov exponent $L_{max}$ of network \dref{n1}, which
can be calculated from system \dref{f3} and is a function of
$\alpha$, is referred to as the \textit{master stability function}.
In addition, the region $S$ of negative real $\alpha$ where
$L_{max}$ is also negative is called the \textit{synchronized
region} of network \dref{n1}. Based on the results of
\cite{lu04,pec98},
 the synchronized solution
 of  network \dref{n1} is asymptotically stable if, and
only if,
\begin{equation}\label{f4}
c\lambda_k\in S,\;k=2,3,\cdots,N.
\end{equation}

 The synchronized region $S$ can be an
unbounded region, a bounded region, an empty set, or a union of
several regions. If the synchronous state is an equilibrium point,
then $Df(s(t))$ and $DH(s(t))$ reduce to constant matrices, denoted
by $F$ and $H$, respectively. In this case, system \dref{f3} becomes
\begin{equation}\label{f5}
\dot{\omega}=[F+\alpha H]\omega.
\end{equation}
 Hence, the synchronized region $S$ becomes the stability region of
$F+\alpha H$ with respect to parameter $\alpha$.  This paper mainly
studies this case when the synchronous state is an equilibrium
point.

The rest of this paper is organized as follows. In Section 2, the
disconnected stability region problem for the matrix pencil
$F+\alpha H$ is studied, where the existence of multiple
disconnected stability regions is theoretically proved. In Section
3, the characteristics of matrix convexity for the stability of
matrix pencils  are discussed, where some conditions for testing the
stability or instability of convex combinations of two vertex
matrices are established. In Section 4, smooth Chua's circuit
networks are simulated to illustrate the theoretical results. The
paper is concluded by the last section.

\section{Disconnected stability regions for matrix pencils}

As discussed in the previous section, when the synchronization state
is an equilibrium state, the synchronized region problem reduces to
a stability problem of the matrix pencil $F+\alpha H$ with respect
to parameter $\alpha$. In this section,  the characteristics of
disconnected stability regions for such matrices are studied. In
order to discuss this problem in the real parameter domain, the
following lemmas are necessary.

{\bf Lemma 1}\,\,If the  real polynomial
 $$p(s)=s^n+\gamma_{n-1}s^{n-1}+\cdots+\gamma_1s+\gamma_0 \quad (\gamma_0>0)$$
  is stable, then for any scalar $\epsilon$,
 $0<\epsilon<\gamma_0,$ the following polynomial
$$p_{\epsilon}(s)=s^n+\gamma_{n-1}s^{n-1}+\cdots+\gamma_1s+\epsilon$$
is stable.

{\bf Proof}\,\,Given that $p(s)$ is stable, polynomial
$p_{\epsilon}(s)$ is stable if, and only if,
 $p(s)-\epsilon$ is stable for all $0<\epsilon<\gamma_0$, or
 equivalently, the function
$$\frac{\frac{\epsilon}{p(s)}}{1-\frac{\epsilon}{p(s)}}$$
is stable. Further, this is equivalent to that the Nyquist plot of
$-\frac{\epsilon}{p(s)}$ does not enclose the point $(-1,0)$ for all
$0<\epsilon<\gamma_0$,
which obviously holds. \hfill $\Box$

{\bf Lemma 2}\,\, Given a polynomial
$p(\alpha)=(\alpha+1)(\alpha+2)\cdots(\alpha+n)$ with variable
$\alpha$ and $n\geq 2$, there is a scalar $\beta>0$ such that
$p(\alpha)-\beta^n$ has $n$ negative real roots.

{\bf Proof}\,\, Take $\beta>0$ such that
$\beta^n<\frac{1}{2}(0.5\times 1.5 \times \cdots \times
([\frac{n}{2}]-0.5))$. Then, one can get
$$p(0)-\beta^n>0,\,p(-1)-\beta^n<0,\,p(-2.5)-\beta^n>0,\,\cdots,\,
p\left(-2\times
\left[\frac{n}{2}\right]-0.5\right)-\beta^n>0,\,p(-n)-\beta^n<0.$$
Therefore, the sign of $p(\alpha)-\beta^n$ changes $n$ times on the
negative real axis. This means that $p(\alpha)-\beta^n$ has $n$ real
roots on the negative real axis.
 \hfill $\Box$

{\bf Lemma 3}\,\, Given two scalars $\beta_0$ and $\beta$ with
$\beta>0$ and $\beta-\beta_0>0$, there are scalars
$0<\alpha_1<\cdots< \alpha_n$ such that
 $\alpha_1\alpha_2\cdots\alpha_n=\beta$ and
 all roots of
$p(\alpha)=(\alpha+\alpha_1)(\alpha+\alpha_2)\cdots(\alpha+\alpha_n)-(\beta-\beta_0)$
are real.

{\bf Proof}\,\, By the method of Lemma 2, it suffices to prove this
lemma by choosing $\alpha_i$ with the above constraints
 such that the sign of $p(\alpha)$
changes $n$ times on the real axis.
 \hfill $\Box$

With the above lemmas, one can get the following  results.

 {\bf Theorem 1}\,\, For any natural number $n$, there are matrices $F$ and $H$ of order
 $2(n-1)$ such that $F+\alpha H$ has $n$ disconnected stable regions with
 respect to  parameter $\alpha$.

{\bf Proof}\,\,  As shown in Lemma 2, one may take $\beta>0$ such
that
 \begin{equation} \label{th1}
  p(\alpha)=(\alpha+1)(\alpha+2)\cdots(\alpha+2(n-1))-\beta^{2(n-1)}=0
  \end{equation}
  has $2(n-1)$ real roots, denoted by $\beta_1, \, \beta_2, \,\cdots, \, \beta_{2(n-1)}.$
  Then, take
  $$H=\left(\begin{array}{cccc}
  0 & -1 & \cdots & 0\\
  \vdots & \ddots & \ddots & \vdots\\
  0&0 &  \ddots  & -1\\
  1 & 0 & \cdots & 0 \end{array}\right),\quad
  F_{1}=\left(\begin{array}{cccc}
  0 & \beta_1 & \cdots & 0\\
  \vdots & \ddots & \ddots & \vdots\\
  0 &0&  \ddots  & \beta_{2n-3}\\
 -\beta_{2(n-1)} & 0 & \cdots & 0 \end{array}\right),
$$
 and $F=-\beta I_{2(n-1)}+F_{1}$, where $I_{2(n-1)}$ is the identity matrix of order
 $2(n-1)$. Obviously, the characteristic
 polynomial of $F+\alpha H$ is
$$\hbox {det} (sI-F-\alpha H)=(s+\beta)^{2(n-1)}+(\alpha-\beta_1)(\alpha-\beta_2)\cdots(\alpha-\beta_{2(n-1)}).$$
 Using (\ref{th1}), one has
 $$\hbox {det} (sI-F-\alpha H)=(s+\beta)^{2(n-1)}-\beta^{2(n-1)}+(\alpha+1)(\alpha+2)\cdots(\alpha+2(n-1)).$$
The constant term in  $\hbox {det} (sI-F-\alpha H)$ is
$$(\alpha+1)(\alpha+2)\cdots(\alpha+2(n-1)), $$
 which is larger than zero if the  parameter $\alpha$ is located in the following $n$
 regions:
 \begin{equation} \label{th12}
 (0, -1),\,(-2, -3),\,  \,\cdots, \, (-2(n-2), -2n+3),\, (-2(n-1), -\infty),
  \end{equation}
 and is smaller than zero if $\alpha$ is located in the following $n-1$
 regions:
 $$(-1, -2),\,(-3, -4),\, \,\cdots, \, (-2n+3, -2(n-1)).$$
 Obviously, by Lemma 1 $\hbox {det} (sI-F-\alpha H)$ has $n$ disconnected stable regions with
 respect to parameter $\alpha$, which are contained in the $n$ regions shown in (\ref{th12}),
 respectively.
 \hfill $\Box$

 Combining with the discussions in Section 1, for any natural number $n$,
 Theorem 1 shows the
 existence of a network which has $n$ disconnected synchronized
 regions. However, for a general network, the node equation is given, i.e., $F$ is given,
 which can not be chosen arbitrarily.
 In this case, one may apply the following result with a chosen
 inner linking matrix $H$.

 {\bf Theorem 2}\,\, For any given real stable matrix $F$ of order
 $n$, suppose $\hbox {det}(sI-F)=s^n+\gamma_{n-1}s^{n-1}+\cdots+\gamma_1
 s+\gamma_0$, and  every eigenvalue of $F$ corresponds to only one Jordan form. If there is a
 scalar $\beta_0\not=0$ such that $p(s)=s^n+\gamma_{n-1}s^{n-1}+\cdots+\gamma_1 s+\gamma_0-\beta_0$
 is stable and $p(s) $ has $n_i$ pairs of conjugate complex eigenvalues, then
  there exists a real matrix $H$ such that $F+\alpha H$ has $[\frac{n-n_i}{2}]+1$
 disconnected stable regions with respect to  parameter $\alpha$.

{\bf Proof}\,\,
 First, suppose that
  there is a scalar $\beta_0$ such that
 $$p(s)=s^n+\gamma_{n-1}s^{n-1}+\cdots+\gamma_1 s+\gamma_0-\beta_0$$
 is  stable, with  $n$ real roots denoted by  $\lambda_{01}, \cdots,
 \lambda_{0n}.$
 Following Lemma 3, take scalars $0<\alpha_1<\cdots< \alpha_n$ such that
 $\alpha_1\alpha_2\cdots\alpha_n=\gamma_0$ and all roots of
 $(\alpha+\alpha_1)(\alpha+\alpha_2)\cdots(\alpha+\alpha_n)-(\gamma_0-\beta_0)$ are real, denoted by
 $-\beta_1, \cdots, -\beta_n$.
 Consequently,
  \begin{equation} \label{th2}
  p(\alpha)=(\alpha+\alpha_1)(\alpha+\alpha_2)\cdots(\alpha+\alpha_n)-(\gamma_0-\beta_0)=
 (\alpha+\beta_1)(\alpha+\beta_2)\cdots(\alpha+\beta_n).
 \end{equation}
 Obviously, $\beta_1\cdots\beta_n=\beta_0.$ Furthermore,
 take
 $$H_{0}=\left(\begin{array}{cccc}
 0 & 1 & \cdots & 0\\
  \vdots & \ddots & \ddots & \vdots\\
  0&0 &  \ddots  & 1\\
  -1 & 0 & \cdots & 0 \end{array}\right),\quad
  F_{0}=\left(\begin{array}{cccc}
  \lambda_{01} & \beta_1 & \cdots & 0\\
  \vdots & \ddots & \ddots & \vdots\\
  0 &0&  \ddots  & \beta_{n-1}\\
 -\beta_{n} & 0 & \cdots & \lambda_{0n}\end{array}\right).
$$
Then,
$\hbox{det}(sI-F_0)=(s-\lambda_{01})\cdots(s-\lambda_{0n})+\beta_0=s^n+\gamma_{n-1}s^{n-1}+\cdots+\gamma_1
 s+\gamma_0=\hbox{det}(sI-F)$.
Hence $F_0$ is similar to $F$, since each eigenvalue of $F$
corresponds to only one Jordan form. Moreover,
$$\hbox{det}(sI-F_{0}-\alpha
H_{0})=(s-\lambda_{01})\cdots(s-\lambda_{0n})+(\alpha+\beta_{1})\cdots(\alpha+\beta_{n})$$
$$=(s-\lambda_{01})\cdots(s-\lambda_{0n})-(\gamma_0-\beta_0)
+(\alpha+\alpha_1)(\alpha+\alpha_2)\cdots(\alpha+\alpha_n).$$ The
constant term in $\hbox{det}(sI-F_{0}-\alpha H_{0})$ is
$$(\alpha+\alpha_1)(\alpha+\alpha_2)\cdots(\alpha+\alpha_n),$$
which is larger than zero if the parameter $\alpha$ is located in
the following $[\frac{n}{2}]+1$
 regions:
 \begin{equation} \label{th21}
 (0, -\alpha_1),\,(-\alpha_2, -\alpha_3),\,  \,\cdots, \,
  \end{equation}
 and is smaller than zero if $\alpha$ is located in the following $n-[\frac{n}{2}]$
 regions:
 $$(-\alpha_1, -\alpha_2),\,(-\alpha_3, -\alpha_4),\, \,\cdots \, .$$
 Thus, by Lemma 1, $F_{0}+\alpha H_{0}$ has $[\frac{n}{2}]+1$ disconnected stable regions, which
 are located in the regions shown in (\ref{th21}).
Since $F_{0}$ is similar to $F$,  there exists a nonsingular matrix
$P$ such that $P^{-1}F_{0}P=F$. Therefore, $H=P^{-1}H_{0}P$ is the
matrix to be found. And $F+\alpha H$ has the same stable regions as
$F_{0}+\alpha H_{0}$.

Then, assume that there are some conjugate complex pairs in
$\lambda_{01}, \cdots, \lambda_{0n}$. For simplicity, suppose that
there is only one pair of conjugate complex eigenvalues,
$\lambda_{01}=\xi_1+\eta_1i, \lambda_{02}=\xi_1-\eta_1i$, and
$\lambda_{03}, \cdots, \lambda_{0n}$ are all real.

Similarly to the above proof, take scalars $0<\alpha_2<\cdots<
\alpha_n$ such that
 $\alpha_2\alpha_2\cdots\alpha_n=\gamma_0$ and all roots of
 $(\alpha+\alpha_2)(\alpha+\alpha_3)\cdots(\alpha+\alpha_n)-(\gamma_0-\beta_0)$ are real,
 denoted by
 $-\beta_2, \cdots, -\beta_n$.
 Consequently,
  \begin{equation} \label{th22}
  p(\alpha)=(\alpha+\alpha_2)(\alpha+\alpha_3)\cdots(\alpha+\alpha_n)-(\gamma_0-\beta_0)=
 (\alpha+\beta_2)(\alpha+\beta_3)\cdots(\alpha+\beta_n).
 \end{equation}
 Obviously, $\beta_2\cdots\beta_n=\beta_0.$
Furthermore, take
 $$H_{0}=\left(\begin{array}{ccccc}
0 & 0 & 0 & \cdots & 0\\
 0 & 0 & 1 & \cdots & 0\\
  \vdots & \vdots & \ddots & \ddots & \vdots\\
  0&0 & 0 & \ddots  & 1\\
  -1 & 0 & 0& \cdots & 0 \end{array}\right),\quad
  F_{0}=\left(\begin{array}{cccccc}
 \xi_1 & 1 & 0 & 0 & \cdots & 0\\
 -\eta_1^2 & \xi_1 & \beta_2 & 0 & \cdots & 0\\
0 & 0 &  \lambda_{03} & \beta_3 & \cdots & 0\\
\vdots& \vdots&  \vdots & \ddots & \ddots & \vdots\\
  0 & 0 & 0 &0&  \ddots  & \beta_{n-1}\\
 -\beta_{n} & 0 & 0 & 0 & \cdots & \lambda_{0n}\end{array}\right).
$$
Obviously,
$\hbox{det}(sI-F_0)=(s-2\xi_1s+\xi_1^2+\eta_1^2)(s-\lambda_{03})\cdots(s-\lambda_{0n})+\beta_0
=s^n+\gamma_{n-1}s^{n-1}+\cdots+\gamma_1
 s+\gamma_0=\hbox{det}(sI-F)$.
Hence, $F_0$ is similar to $F$. Moreover,
$$\hbox{det}(sI-F_{0}-\alpha
H_{0})=(s-2\xi_1s+\xi_1^2+\eta_1^2)(s-\lambda_{03})\cdots(s-\lambda_{0n})+(\alpha+\beta_{2})\cdots(\alpha+\beta_{n})$$
$$=s^n+\gamma_{n-1}s^{n-1}+\cdots+\gamma_1 s
+(\alpha+\alpha_2)\cdots(\alpha+\alpha_n).$$ The constant term in
$\hbox{det}(sI-F_{0}-\alpha H_{0})$ is
$$(\alpha+\alpha_2)\cdots(\alpha+\alpha_n),$$
which is larger than zero if the parameter $\alpha$ is located  in
the following $[\frac{n-1}{2}]+1$
 regions:
 \begin{equation} \label{th23}
 (0, -\alpha_2),\, (-\alpha_3, -\alpha_4), \,\cdots.
  \end{equation}
Repeating the process as above, one can complete the proof easilly.
 \hfill $\Box$

 {\bf Remark 1}\,\, For simplicity, in both Theorems 1 and 2,  the parameter
 $\alpha$ only appears in the constant term of the characteristic
 polynomial of $F+\alpha H$. In this way,  the constant term in
$ \hbox{det}(sI-F-\alpha H)$, which is a polynomial of parameter
$\alpha$, simply determines the disconnected stable regions. If
$\alpha$ appears in the higher-order terms of $
\hbox{det}(sI-F-\alpha H)$, the problem becomes harder to solve,
leaving  an interesting  topic for further research.

{\bf Remark 2}\,\, In order to guarantee  $H$ be a real matrix, two
cases are considered in the proof of Theorem 2, i.e., there are or
there are no conjugate complex  pairs in $\lambda_{0i},
i=1,\cdots,n$. If $H$ can be chosen to be  a complex matrix, the
proof of Theorem 2 will be simplified and  $H$ may be chosen such
that $F+\alpha H$ has $[\frac{n}{2}]+1$ disconnected stable regions.
In addition,  if all $\lambda_{0i}, i=1,\cdots,n$, are complex
scalars, then there exists a real $H$ such that $F+\alpha H$ has at
least $[\frac{n}{4}]+1$
 disconnected stable regions.

 According to the above discussions, one can also choose a suitable $H$
 such that $F+\alpha H$ has only one convex stable region with
 respect to parameter $\alpha$, as further discussed below.

\section{Characteristics of convexity for stability of matrix pencils }

In the previous section, it shows the existence of any $n$
disconnected stable regions of the matrix pencil $F+\alpha H$.
Contrary to this non-convexity, given two parameter values
$\alpha_1$ and $\alpha_2$, whether or not the stability of
$F+\alpha_1 H$ and $F+\alpha_2 H$ implies the stability of
$F+(\lambda\alpha_1+(1-\lambda)\alpha_2) H$, for all $0\leq
\lambda\leq 1$, is  an interesting problem. Obviously, a good
understanding of this convexity characteristic is useful for
enhancing the stability of the matrix pencil  $F+\alpha H$.

{\bf Lemma 4}\,\, Suppose that $F+\alpha_1 H$ and $F+\alpha_2 H$ are
stable, and the rank of $H$ is 1. Let $H=bc$, where $b$ is a column
vector and $c$ is a row vector with compatible dimensions, and $(F,
b)$ be controllable. Then the following conditions are equivalent to
each other:

(i) $\lambda (F+ \alpha_1 H)^{-1}+ (1-\lambda) (F+\alpha_2 H)$ is
stable for all $0\leq\lambda \leq 1$.

(ii) There is a common matrix $P=P^T$ such that
$$P(F+\alpha_i H)+(F+\alpha_i H)^TP<0, \quad i=1, 2.$$

(iii) $(F+\alpha_1 H)(F+\alpha_2 H)$ does not have negative real
eigenvalues.

(iv) $1-{\bf Re}\{(\alpha_2-\alpha_1)c(jw I-F-\alpha_1
H)^{-1}b\}>0$, $\forall \, w\in {\bf R}$. \\
Further, if any one of (i)-(iv) holds, one has

($\star$) $\lambda(F+ \alpha_1 H)+ (1-\lambda)(F+ \alpha_2 H)$ is
stable for all $0\leq\lambda \leq 1.$


{\bf Proof} \,\,See \cite{sho04, sho03} for the equivalences among
(ii)-(iv). Now, if (i) holds, then
 $$\hbox {det}(\lambda (F+ \alpha_1 H)^{-1}+ (1-\lambda) (F+\alpha_2
 H))\not=0,\,\, \forall \,\,0\leq\lambda \leq 1.$$
Equivalently,
$$\hbox {det}\left( \frac{\lambda}{1-\lambda} I+ (F+ \alpha_1 H) (F+\alpha_2
 H)\right)\not=0,\,\, \forall \,\,0\leq\lambda \leq 1, $$
which implies (iii). On the other hand, by a simple congruence
transformation with matrix $(F+\alpha_1 H)^{-1}$,  (ii) implies the
existence of a common matrix $P=P^T$ such that
 \begin{equation} \label{le4}
P(F+\alpha_1 H)^{-1}+(F+\alpha_1 H)^{-T}P<0, \,\,P(F+\alpha_2
H)+(F+\alpha_2 H)^{T}P<0,
 \end{equation} which implies (i).  And,
obviously, (ii) implies ($\star$). This completes the proof. \hfill
$\Box$

{\bf Remark 3}\,\, Any one of Lemma 4 (i)-(iv) implies ($\star$).
However, generally, ($\star$) does not imply the other conditions of
Lemma 4. For example, with
$$F=\left(\begin{array}{cc}0 & 1\\-1 &
-0.1\end{array}\right),\,
 H=\left(\begin{array}{cc}0 & 0\\1 & 0\end{array}\right),\,
 \alpha_1=0, \alpha_2=0.9,$$
 obviously Lemma 4 ($\star$) holds for the above matrices, but
 all Lemma 4 (i)-(iv) do not hold.

{\bf Remark 4}\,\,Obviously, Lemma 4 (ii) is equivalent to
(\ref{le4}).  However, since the rank of $(F+\alpha_1
H)^{-1}-(F+\alpha_2 H)$ is generally not 1, (\ref{le4}) is not
equivalent to that $(F+\alpha_1 H)^{-1}(F+\alpha_2 H)$ does not have
negative real eigenvalues. In fact, $(F+\alpha_1 H)^{-1}(F+\alpha_2
H)=(F+\alpha_1 H)^{-1}(F+\alpha_1 H+(\alpha_2-\alpha_1)
H)=I+(\alpha_2-\alpha_1)(F+\alpha_1 H)^{-1} H$.  That  $(F+\alpha_1
H)^{-1}(F+\alpha_2 H)$ does not have negative real eigenvalues is
equivalent to $1+(\alpha_2-\alpha_1)c(F+\alpha_1 H)^{-1}b>0$, which
is Lemma 4 (iv) at $jw=0$. The example given in Remark 3 can also be
an example for this.

{\bf Remark 5}\,\, It is not necessary to require $P>0$ in Lemma 4
(ii). The positive definiteness of $P$ is naturally guaranteed by
the stability of $F+\alpha_i H, i=1,2.$  For this property, the
following lemma is useful.

{\bf Lemma 5}\,\cite{leo96}\, Suppose that $F$ and $P=P^T$ are
matrices of order $n$. If $PF+F^TP<0$, then $F$ has no eigenvalues
on the imaginary axis, det$(P)\not=0$, and the number of positive
eigenvalues of $P$ is equal to the number of eigenvalues of $F$ with
negative real parts.

 {\bf Corollary 1}\,\, Suppose the rank of $H$ is 1 and
$F+\alpha_1 H$ is stable. Let $H=bc$ and $(F,b)$ be controllable.
Then, $F+\alpha H$ is stable for all $\alpha\in(-\infty, \alpha_1]$
if, and only if, ${\bf Re}\{c(jw I-F-\alpha_1 H)^{-1}b\}\leq 0,$
$\forall \, w\in {\bf R}.$

As to the instability of matrix pencils, the following result can be
obtained from Lemmas 4 and  5.

 {\bf Theorem 3}\,\, Suppose that $F+\alpha_1 H$ and $F+\alpha_2 H$ are unstable and do not have
 imaginary eigenvalues, and the rank of $H$ is 1. Let $H=bc$ and $(F, b)$ be controllable. Then,
 the conditions of  Lemma 4 (ii)-(iv) are equivalent to each other. Further, if
 any one of Lemma 4 (ii)-(iv) holds, then $F+(1-\lambda)\alpha_1H+\lambda \alpha_2H$ is
unstable for all $\lambda\in (0, 1).$

{\bf Proof}\,\, Since the Kalman-Yakubovich-Popov lemma also holds
for unstable state matrices \cite{ran96},  Theorem 5 in
\cite{duan05} yields the equivalence between (ii) and (iv) of Lemma
4. Since $F+\alpha_1 H$ and $F+\alpha_2 H$  do not have
 imaginary eigenvalues, Theorem 3.1 in \cite{sho03} gives the
 equivalence between (ii) and (iii) of Lemma 4.

In what follows, it is to prove that Lemma 4 (iv) implies the
instability of $F+\alpha_1H+\lambda (\alpha_2-\alpha_1)H$.
 Suppose that there is a $\lambda_0 \in (0, 1)$ such
that $F+\alpha_1H+\lambda_0 (\alpha_2-\alpha_1)H$ is stable. Since
matrix eigenvalues change continuously with matrix parameters, by
the instability of $F+\alpha_1 H$, there exists a $\lambda_1, \,
0<\lambda_1<\lambda_0$ such that $F+\alpha_1H+\lambda_1
(\alpha_2-\alpha_1)H$ has an imaginary eigenvalue $jw_0$, i.e.,
$$\hbox{det}(jw_0I-F-\alpha_1H-\lambda_1
(\alpha_2-\alpha_1)H)=0.$$
 Since $F+\alpha_1 H$  does not
have imaginary eigenvalues, the above condition is equivalent to
$$\hbox{det}\left( \frac{1}{\lambda_1}I-
(\alpha_2-\alpha_1)(jw_0I-F-\alpha_1H)^{-1}H\right)=0,$$ or
$$ \frac{1}{\lambda_1}-
(\alpha_2-\alpha_1)c(jw_0I-F-\alpha_1H)^{-1}b=0,$$ which is contrary
to Lemma 4 (iv). This completes the proof.  \hfill $\Box$

{\bf Remark 6}\,\, The common Lyapunov matrix problem for stable
matrix pencils was studied in \cite{ sho04, sho03}. Theorem 3 above
generalizes the similar results to unstable matrix pencils. In fact,
any one of Lemma 4 (ii)-(vi) guarantees that transferring an
eigenvalue  between the left-half and  right-half complex planes is
impossible. Therefore, when any one of Lemma 4 (ii)-(vi) holds,
$F+\alpha_1 H$ and $F+\alpha_2 H$ must have the same number of
eigenvalues with positive real parts. In addition, the stability of
the convex combinations of $(F+\alpha_1 H)^{-1}$ and $F+\alpha_2 H$
is equivalent to any one of Lemma 4 (ii)-(iv). However, for
instability, this equivalence does not hold; that is, the
instability of $(F+\alpha_1H)^{-1}+\lambda (\alpha_2-\alpha_1)H,$
$\lambda\in [0, 1],$ does not necessarily imply any one of Lemma 4
(ii)-(vi). For example, with
$$F=\left(\begin{array}{cc}1 & 1\\0 &
 1\end{array}\right),\,
 H=\left(\begin{array}{cc}0 & 0\\0 & 1\end{array}\right),\,
 \alpha_1=0, \alpha_2=-2,$$
   obviously, $(F+\alpha_1H)^{-1}+\lambda (\alpha_2-\alpha_1)H$ is
  unstable for all $\lambda\in [0, 1],$ but any one of Lemma 4
  (ii)-(iv) does not hold.

Although the characteristics of convexity for stability or
instability of the matrix pencil $F+\alpha H$ have been discussed
when the rank of $H$ is 1, it is still hard  to decide the stability
or instability of the convex combinations of $F+\alpha_1 H$ and
$F+\alpha_2 H$ for a general $H$. The results of Lemma 4 and Theorem
3 are directly related to the existence of a common matrix $P$ for
two vertex matrices. For a general $H$, this common-matrix  method
is very conservative. In this case, nevertheless, the following
lemma provides a less conservative criterion \cite{pea00}.

{\bf Lemma 6}\,\,Suppose that $F+\alpha_1 H$ and $F+\alpha_2 H$ are
stable. If there are matrices $P_1=P_1^T, P_2=P_2^T$, $G$ and  $V$
such that
$$\left(\begin{array}{cc}-G-G^T & P_i-V^T+G(F+\alpha_i H)\\
P_i-V+(F+\alpha_i H)^TG^T & V(F+\alpha_i H)+(F+\alpha_i H)^TV^T
\end{array}\right)<0, \,\, i=1,2, $$ then
 $F+\lambda \alpha_1 H+ (1-\lambda) \alpha_2 H$ is stable
for all $0\leq\lambda \leq 1$. \hfill $\Box$

In the above lemma, by introducing new slack matrices $G$ and $V$,
the symmetrical matrices $P_1$ and $P_2$ can be chosen
parameter-dependent for the study of stability of the convex
combination of $F+\alpha_1 H$ and $F+\alpha_2 H$. Similarly to the
method used in \cite{duan06, liu07}, one can also design controllers
to enlarge stability regions by the above-discussed  method.

For instability of matrix pencils, one can likewise obtain the
following result.

{\bf Theorem 4}\,\,Suppose that $F+\alpha_1 H$ and $F+\alpha_2 H$
are unstable.
 If there are matrices $P_1=P_1^T,
P_2=P_2^T$, $G$ and $ V$ such that
$$\left(\begin{array}{cc}-G-G^T & P_i-V^T+G(F+\alpha_i H)\\
P_i-V+(F+\alpha_i H)^TG^T & V(F+\alpha_i H)+(F+\alpha_i H)^TV^T
\end{array}\right)<0, \,\, i=1,2, $$ then
 $F+(\lambda \alpha_1 + (1-\lambda) \alpha_2) H$ is unstable
for all $0\leq\lambda \leq 1$.

{\bf Proof }\,\,For any $\lambda \in(0, 1)$, taking a convex
combination between two inequalities in Theorem 4, one gets
 \begin{equation} \label {th41}
  \left(\begin{array}{cc}-G-G^T & \lambda P_1+(1-\lambda) P_2-V^T+GF_{\lambda}\\
 \lambda P_1+(1-\lambda) P_2-V+F_{\lambda}^TG^T & VF_{\lambda}+F_{\lambda}^TV^T
\end{array}\right)<0, \end{equation}
where $F_{\lambda}=F+(\lambda\alpha_1+(1-\lambda)\alpha_2) H.$ Then,
by the similar method used in \cite{pea00}, the above inequality is
equivalent to
$$(\lambda P_1+(1-\lambda) P_2)F_{\lambda}+F_{\lambda}^T(\lambda P_1+(1-\lambda)
P_2)<0.$$
 Here,  $F_{\lambda}$ must be unstable. If it was stable, then, as
proved in
 Theorem 3, there would exist $\lambda_0\in (0, \lambda)$ such that
  $F_{\lambda_0}=F+(\lambda_0\alpha_1+(1-\lambda_0)\alpha_2) H$
  has an imaginary eigenvalue $jw_0$. Then, by Lemma 5,
$$(\lambda_0 P_1+(1-\lambda_0)
P_2)F_{\lambda_0}+F_{\lambda_0}^T(\lambda_0 P_1+(1-\lambda_0) P_2)$$
can not be strictly negative definite, which is contrary to
(\ref{th41}).
  \hfill $\Box$

{\bf Remark 7}\,\, Theorem 4 generalizes the method of \cite{pea00}
to the instability of matrix pencils. Obviously, the instability of
the matrix pencil $F+\alpha H$ is important  in desynchronization
problems.
  As discussed in Remark 6, if
Theorem 4 holds, transferring any eigenvalue of $F+\alpha_1 H$ and
$F+\alpha_2 H$ between the left-half and right-half complex planes
is impossible. Therefore, $F+\alpha_1 H$ and $F+\alpha_2 H$ must
have the same number of eigenvalues with positive real parts when
Theorem 4 holds.

\section{Synchronization of  smooth Chua's circuit networks}

In this section, consider the synchronization problem in a network
of  smooth Chua's circuits.

{\bf Example 1}\,\, Consider the network (\ref{n1}) consisting of
the third-order smooth Chua's circuits \cite{tsu05}, in which each
node equation is
 \begin{equation}  \label {sc1}
\begin{array}{ccl}
 \dot{x}_{i1}&=&-k\alpha x_{i1}+k\alpha x_{i2}-k\alpha(a x_{i1}^3+bx_{i1}),\\
 \dot{x}_{i2}& =& kx_{i1}-kx_{i2}+kx_{i3},\\
 \dot{x}_{i3}& =& -k\beta x_{i2}-k\gamma x_{i3}. \end{array}\end{equation}
The vector  $x_i$ in (\ref{n1}) is $(x_{i1}, x_{i2}, x_{i3})^T$
here. Linearizing (\ref{sc1})  at its zero equilibrium gives
 \begin{equation}  \label {lsc1}
 \dot{x}_{i}= Fx_i, \quad F=\left(\begin{array}{ccc} -k\alpha-k\alpha b &k\alpha & 0\\
 k & -k & k\\0 & -k\beta &-k\gamma  \end{array}\right).
 \end{equation}

 Take $k=1,
\alpha=-0.1, \beta=-1, \gamma=1,  a=1, b=-25.$ Then $F$ is stable,
i.e., the node system (\ref{sc1}) is locally stable about zero. One
can easily take a parameter $\beta_0=-0.8$ such that all roots of
$\hbox{det}(sI-F)-\beta_0$ are real. Following the method of Theorem
2, take $\alpha_1=0.01,  \alpha_2=1, \alpha_3=10,$ and
 $$H=\left(\begin{array}{ccc} 0.8348 &  9.6619&  2.6591\\ 0.1002 & 0.0694 & 0.1005\\
  -0.3254 &  -8.5837 &  -0.9042  \end{array}\right).$$
  Then, by simply computation, one knows that $F+\alpha H$ has two
  disconnected stable regions: $S_1=[-0.0099, 0]$ and $S_2=[-2.225, -1)$.
  Therefore, the entire synchronized region   is $S_1\bigcup S_2$. Further,
   suppose that the number of nodes is $N=10$, and the outer coupled matrix $A$ is a
  globally coupled matrix, i.e., all the diagonal entries of $A$ are
  $-9$ and the other entries are 1, which has eigenvalues
  $$\lambda_1=0, \lambda_2=\cdots=\lambda_{10}=-10.$$
  Then, by (\ref{f4}), network (\ref{n1}) with the above parameter values achieves
  local synchronization when the coupling strength $c$ satisfies $c\in
  [0, 0.00099]$ or  $c\in(0.1, 0.2225]$. Figures 1 and 2
  show the synchronization and non-synchronization behaviors of this network.

\begin{center}
 \unitlength=1cm
 \qquad \hbox{\hspace*{0.1cm}  \epsfxsize5cm \epsfysize5cm
\epsffile{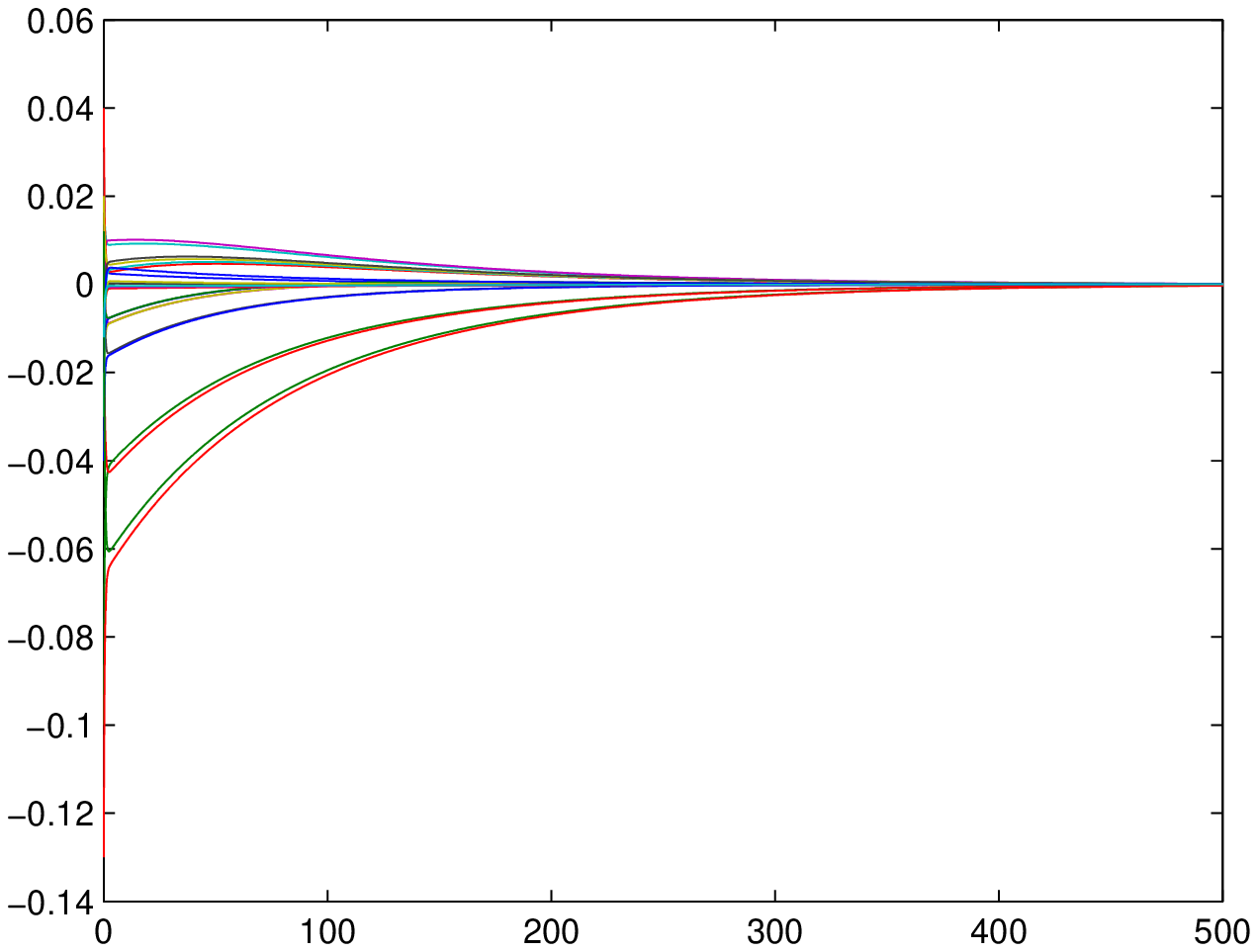}\qquad \epsfxsize5cm \epsfysize5cm
\epsffile{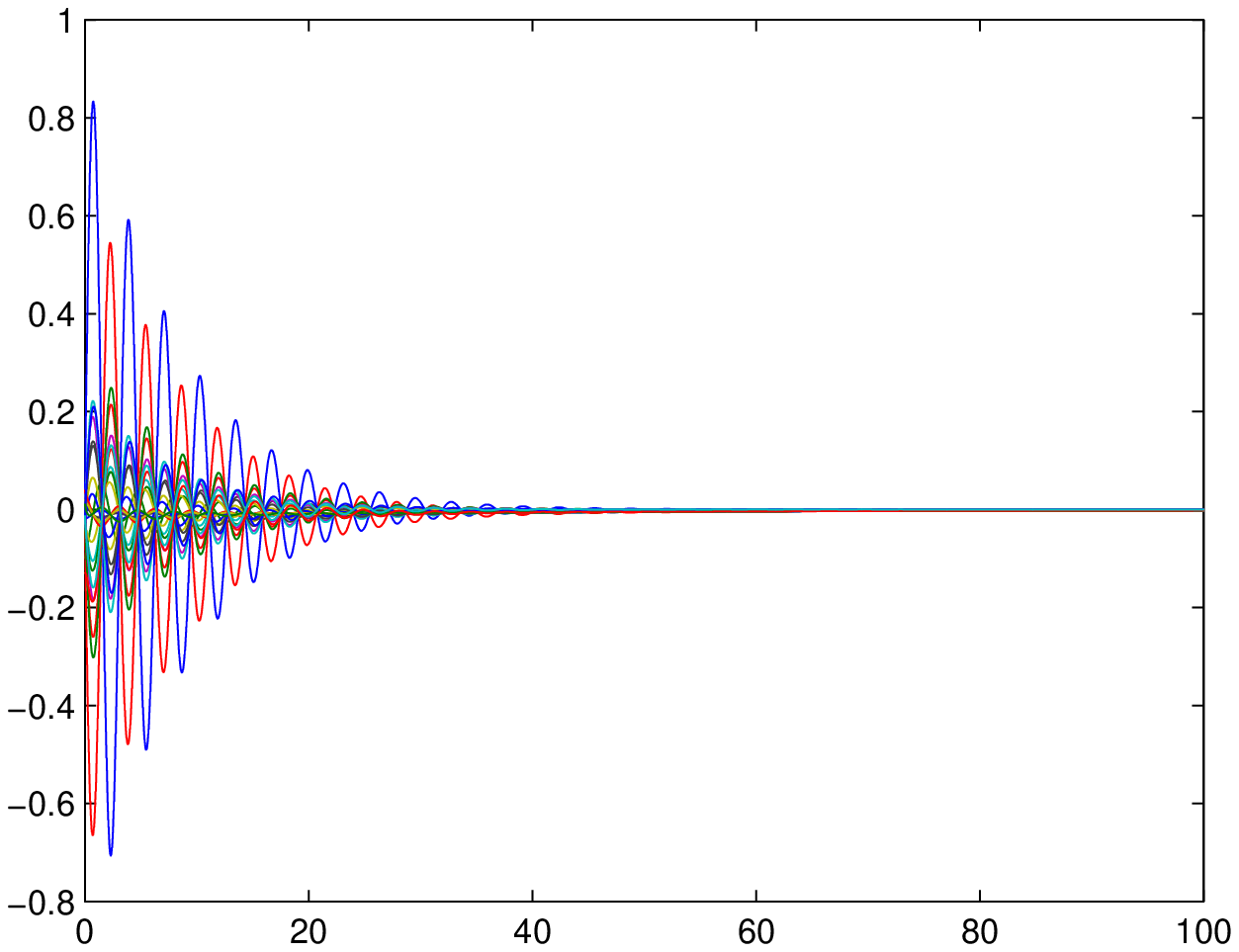}
 }
\end{center}
\vskip -0.3cm \centerline{ (a) $c=0.0005\in[0, 0.00099].$ \qquad
\quad(b) $c=0.2\in(0.1, 0.2225]$.}
 \centerline{ Fig. 1 \,\, Synchronization of network (\ref{n1}) with different coupling strengths.}

\begin{center}
 \unitlength=1cm
 \qquad \hbox{\hspace*{0.1cm}  \epsfxsize5cm \epsfysize5cm
\epsffile{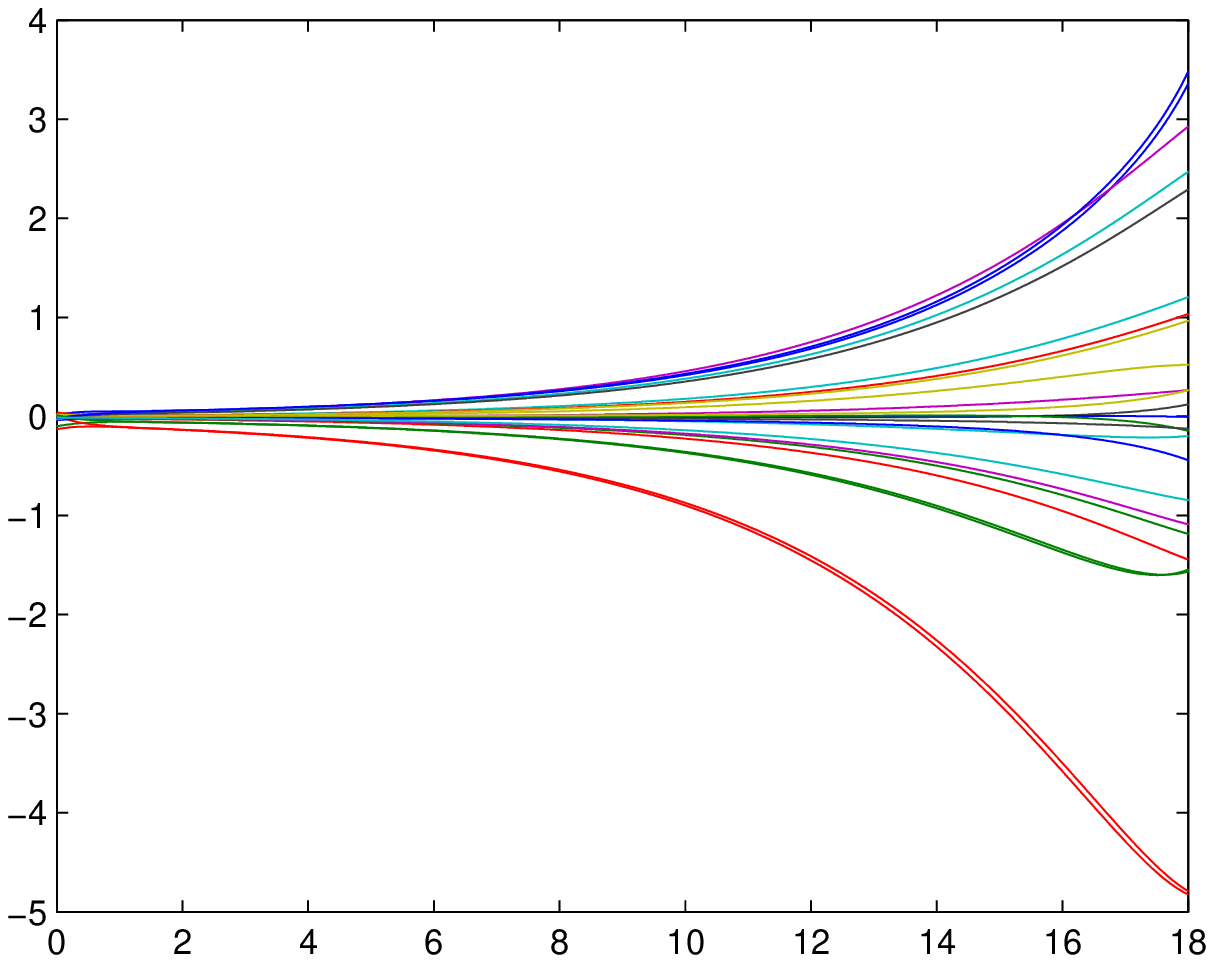}\qquad \epsfxsize5cm \epsfysize5cm
\epsffile{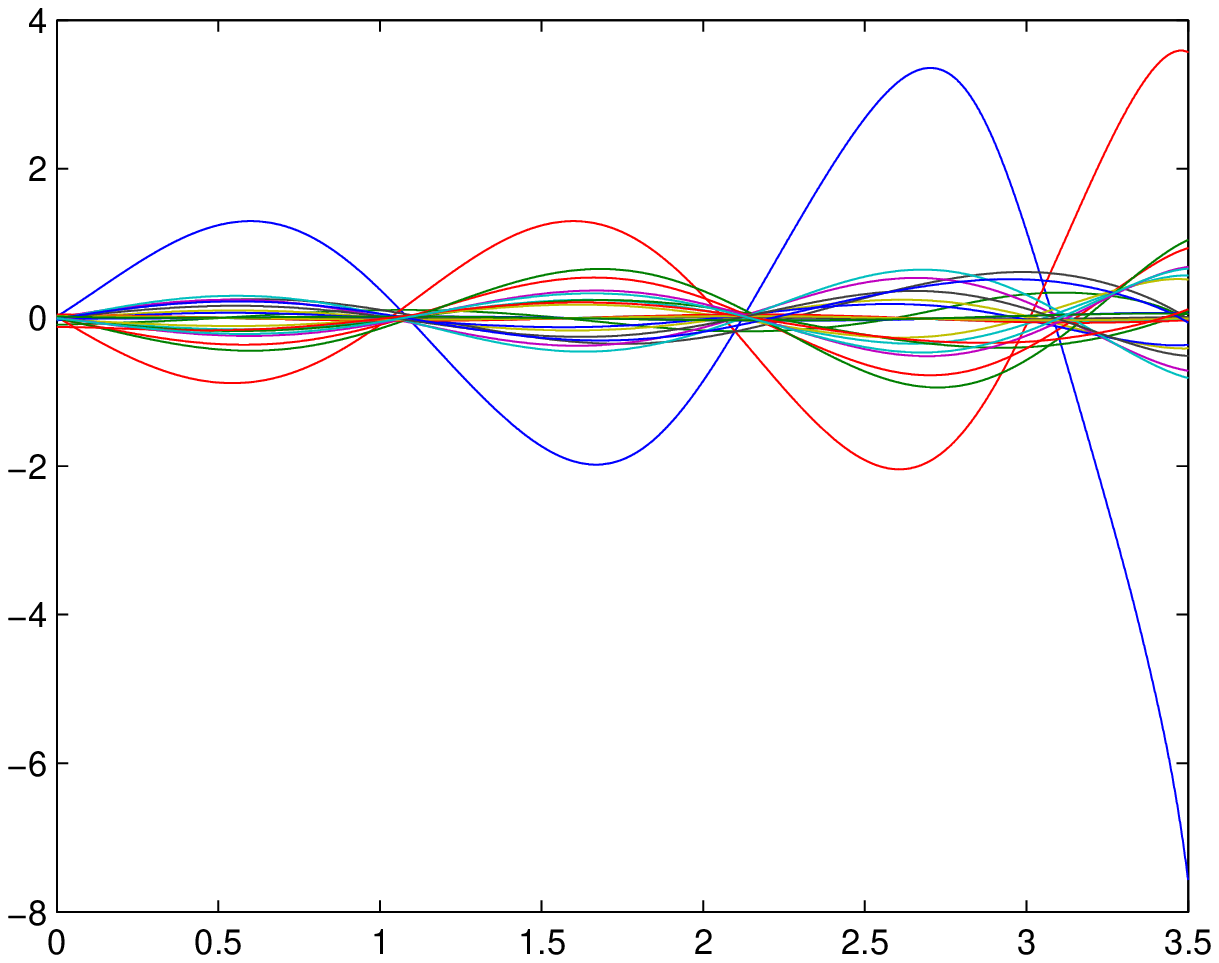}
 }
\end{center}
\vskip -0.3cm \centerline{ (a) $c=0.02\in(0.001, 0.1).$ \qquad
\quad(b) $c=0.3\in(0.2225, +\infty)$.}
 \centerline{ Fig. 2 \,\, Non-synchronization of network (\ref{n1}) with different coupling strengths.}

{\bf Example 2}\,\, Consider the network (\ref{n1}) consisting of
the forth-order generalized smooth Chua's circuits \cite{duan07},
with node equation
 \begin{equation}  \label {sc2}
\begin{array}{ccl}
 \dot{x}_{i1}&=&-k\alpha x_{i1}+k\alpha x_{i2}-k\alpha(a x_{i1}^3+bx_{i1}),\\
 \dot{x}_{i2}& =& kx_{i1}-kx_{i2}+kx_{i3},\\
 \dot{x}_{i3}& =& k\beta x_{i2}+k\gamma x_{i4},\\
 \dot{x}_{i4}& =& -0.1 x_{i2}. \end{array}\end{equation}
The vector $x_i$ in (\ref{n1}) is $(x_{i1}, x_{i2}, x_{i3},
x_{i4})^T$ here. Linearizing (\ref{sc2})  at its zero equilibrium
yields
 \begin{equation}  \label {lsc2}
 \dot{x}_{i}= Fx_i, \quad F=\left(\begin{array}{cccc} -k\alpha-k\alpha b &k\alpha & 0 &0\\
 k & -k & k & 0\\0 & k\beta & 0&k\gamma\\0 & -0.1 & 0 & 0  \end{array}\right).
 \end{equation}

 Take parameters $k=3,
\alpha-0.1, \beta=-0.2, \gamma=0.2,  a=1$ and $b=-25$. Then $F$ in
(\ref{lsc2}) is stable, i.e., the node system (\ref{sc2}) is locally
stable about zero. One can easily take a parameter $\beta_0=-0.1$
such that all roots of $\hbox{det}(sI-F)-\beta_0$ are real.
Following  the method of Theorem 2, take $\alpha_1=0.1,
\alpha_2=0.5, \alpha_3=2, \alpha_4=12.96,$ and
 $$H=\left(\begin{array}{cccc} 0.8442 & 0.6319&  0.3547&  -1.8905\\-12.7738 & -9.9676 &19.4669& -20.2986\\
 -10.3570& -8.3421& 18.4474& -20.5913\\-4.7028 &-3.8156& 8.5403& -9.3240  \end{array}\right).$$
  Then, by simply computation, one knows that $F+\alpha H$ has three
  disconnected stable regions: $S_1=(-0.1, 0]$, $S_2=(-2, -0.5)$ and
  $S_3=[-12.95, -12.94]$,
  so the whole synchronized region
  is $S_1\bigcup S_2\bigcup S_3$. Theorem 2 implies that, the
  region $S_3$ should be contained in $(-\infty, -12.96),$ but
  due to the  computing error $S_3$ becomes $[-12.95, -12.94]$, slightly off-set from the
  theoretical prediction.

  Similarly to Example 1, if  $N=10$ and  the outer coupled matrix $A$ is a
  globally coupled matrix, then both  synchronization
  and non-synchronization phenomena can be discussed.

\section{Conclusion}

In this paper, the problem of disconnected synchronized regions has
been carefully studied. When the synchronization state is an
equilibrium point, the problem is reduced to the stability problem
of matrix pencils. The existence of multiple disconnected
synchronized regions is theoretically proved for  network with
higher-dimensional nodes. Further, the characteristics of convexity
for matrix pencils has been discussed. Some test conditions for
stability and instability of convex combinations of vertex matrices
have also been established. Finally, networks of smooth and
generalized smooth Chua's circuits have been simulated to illustrate
the analytic results.

 \hspace*{34pt}

\end{document}